\documentclass[12pt]{amsart}

\newtheorem{theorem}{Theorem}[section]
\newtheorem{lemma}[theorem]{Lemma}

\usepackage[dvips]{graphics}

\theoremstyle{definition}

\theoremstyle{remark}

\theoremstyle{conjecture}

\theoremstyle{corollary}

\theoremstyle{problem}

\numberwithin{equation}{section}

\large\normalsize

\begin{document}

\title{the enumeration of independent sets on some lattices}

\author{Zuhe Zhang}
\address{ Department of Mathematics and Statistics,
The University of Melbourne,
Parkville, VIC 3010, Australia}
\email{Zhang.Zuhe@gmail.com}

%; Secondary 46E25, 20C20}


\begin{abstract}
In this paper, firstly we show that the entropy constants of the number of independent sets on certain plane lattices
are the same as the entropy constants of the corresponding cylindrical and toroidal lattices. Secondly,
we consider three more complex lattices which can not be handled by a single transfer matrix as in the plane quadratic lattice case.
By introducing the concept of transfer multiplicity, we obtain the lower and upper bounds of the entropy constants
of crossed quadratic lattice, generalized aztec diamond lattice and 8-8-4 lattice.
\end{abstract}

\maketitle

\section{Introduction}
The study of lattice statistics in statistical physics has a long history. A
typical problem is to count the ways of putting particles in the sites of a
plane lattice such that no two share the same site or are in
adjacent sites. Such problems are called the planar lattice gases
models \cite{b,b2,bet,dg,f,ps}. Mathematicians formulated them by
the enumeration of the $(0,1)$ matrices which describe the
independent sets in a plane quadratic lattice graph (also called a
planar grid graph) \cite{cw,e,w,w2}. Let us recall some basic concept
of lattice gases model on plane quadratic lattice graph.
We use $G_{m,n}$ (where $m,n$ always denote positive integers) to denote a
plane lattice graph, that is a finite part of a plane lattice, whose
vertices are arranged in $(m+1)$ rows and $(n+1)$ columns. Given an
independent set $S$ of  graph $G_{m,n}$,
 a
portion of $S$ that lies in a fixed column of $G_{m,n}$ can be
represented by an $(m+1)$-vector of $0$'s and $1$'s,  where a $1$
indicates the vertex is in $S$ and a $0$ indicates the vertex is not
in $S$.

If $G_{m,n}$ is a plane quadratic lattice graph, then any
$(m+1)$-vector arising  this way
 has no two consecutive
 $1$'s.
 Clearly, a vertex subset $S$ of $G_{m,n}$ is an
independent set only if all its $n+1$ corresponding $(m+1)$-vectors
represents independent sets. Let $P_{m}$ denote the set of all
$(m+1)$-vectors of $0$'s and $1$'s with  no two consecutive
 $1$'s. The cardinality of $P_{m}$ is well known (and can be easily seen) to be $F_{m+3}$,
the Fibonacci number (starting from $F_0=0, F_1=1$.)
 We can construct an $(m+1)$ by $(n+1)$ matrix $M$ to represent an independent set $S$ of
 $G_{m,n}$
 by   the gluing procedure described below.
 First, take a vector from $P_{m}$ such that it corresponds to the first column of $S$,
 and denote it as  $v_1$. Then, to the right of $v_1$ we glue a  vector $v_2$ selected from
 $P_{m}$
 that corresponds to the second column of $S$. Then we  glue a  vector $v_3$ from $P_{m}$ to
 the right
 of  $v_2$, and so on so forth. Continue this way until the $(n+1)$-th column is glued and
 then we obtain the $(m+1)$ by $(n+1)$ matrix $M$
  representing  the independent set $S$.
 We can get all the possible matrices representing
independent sets of $G_{m,n}$
  by this  procedure of gluing vectors from $P_{m}$. Note that two vectors  of
 $P_{m}$ can be glued together if and only if  they have no 1's in common position,
 i.e., they are orthogonal (their dot product equals zero.)
  In the above procedure, we glue columns from left to right. Similarly, we can have another
   procedure that   glues rows from top to bottom.
  Note that for non-grid lattice graphs $G_{m,n}$,  these two procedures may lead to
  different transfer  patterns.  But for all the lattices considered in this paper,
  the transfer patterns in the two procedures are the same.

\begin{figure}[htbp]
  \centering
 \scalebox{0.4}{\includegraphics{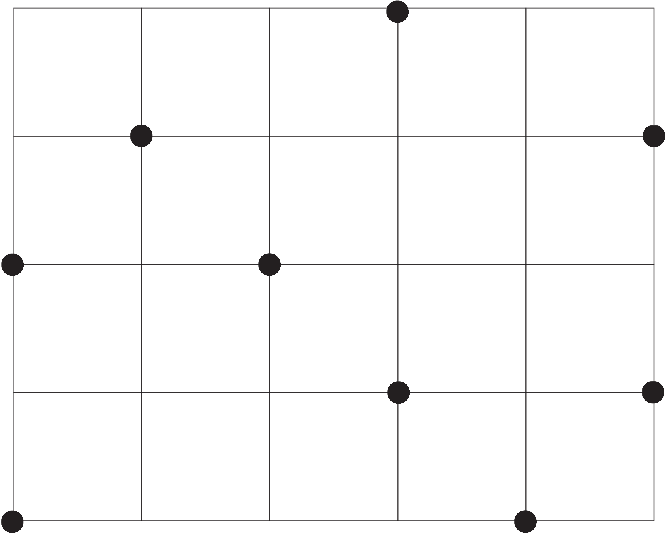}}
  \caption{}\label{fig:fig1.1}
\end{figure}

Fig.\ref{fig:fig1.1} shows an independent set $S$ in the plane
quadratic lattice graph $G_{4,5}$. The portions of $S$ that lie in
each of the columns can be represented by the respective $5$-vectors
$(0,0,1,0,1)$,$(0,1,0,0,0)$,$(0,0,1,0,0)$,$(1,0,0,1,0)$,$(0,0,0,0,1)$,$(0,1,0,1,0)$.

  %%%%%%%%%%%%%%%%%%%%%%%%%%%%%%%%%%%%%%%%%%
%%%%%%%%%%%%%
%%%%%%%%%%%%%%%%%%%%%%%%%%%%%%%%%%%%%%%%%%

%  \begin{figure}[h]
%\center
%  % Requires \usepackage{graphicx}
 % \includegraphics[scale=0.5]{./eps/fig1.1.eps}\\
 % \caption{}\label{fig:fig1.1}
%\end{figure}

  From a plane quadratic
lattice graph $G_{m,n}$ ($n>1$),  identifying its  $m$ edges on the
left with its  $m$ edges on the right, correspondingly, we get a
cylindrical quadratic lattice graph $H_{n,m}$ (Note that here the
graph can be seen as drown on a vertical cylinder.)
 Let $C_{n}$ denote  the set of
$n$-vectors of 0's and 1's  with the property that no two
consecutive $1$'s occur in cyclic order.   Then, similar to the
discussion on $G_{m,n}$, we can see that any independent set $S$ of
$H_{n,m}$ can be represented by an  $(m+1)$ by $n$ matrix $M$ whose
rows are from the
 set $C_{n}$, and that all the representing matrices $M$ can be obtained by the row gluing
 procedure (from top to bottom) similar to the column gluing procedure (from left to right)
 for  $G_{m,n}$. It is not difficult to see that the cardinality of
 the set $C_{n}$ is $F_{n-1}+F_{n+1}$.

 Similarly, by identifying the top and bottom boundary cycles of a
 cylindrical quadratic lattice graph $H_{n,m}$, we  get  a toroidal quadratic lattice graph
 $S_{n,m}$
 that can be seen as drawn on a torus.

In general, for any plane lattice graph $G_{m,n}$, the cylindrical
lattice graph and the toridal lattice graph obtained by the
identifications as above will also be denoted  $H_{n,m}$ and
$S_{n,m}$.

For a given plane quadratic lattice graph $G_{m,n}$,  the transfer
matrix $T_{m}$ is
 an $F_{m+3}\times F_{m+3}$ matrix of $0$'s and
$1$'s, defined as follows.  The rows and columns of $T_{m}$ are
indexed by vectors of $P_{m}$, and the entry of $T_{m}$ in position
$(\alpha$,$\beta)$ is $1$ if the vectors $\alpha$,$\beta$ are
orthogonal, and is $0$ otherwise. Note that the matrix depends only
on $m$.

For example when $m=3$, the possible column vectors of $P_{m}$ are
$(0,0,0,0)$, $(1,0,0,0)$, $(0,1,0,0)$, $(0,0,1,0)$, $(0,0,0,1)$,
$(1,0,1,0)$, $(1,0,0,1)$, $(0,1,0,1)$. If we index the rows and
columns in this order, then the transfer matrix of $G_{m,n}$ is

\begin{center}
$T_{3}=\left(\begin{array}{cccccccc}
1 & 1 & 1 & 1 & 1 & 1 & 1 & 1\\
1 & 0 & 1 & 1 & 0 & 1 & 0 & 1\\
1 & 1 & 0 & 1 & 1 & 1 & 1 & 0\\
1 & 1 & 1 & 0 & 0 & 1 & 1 & 1\\
1 & 0 & 1 & 0 & 0 & 1 & 0 & 1\\
1 & 1 & 1 & 1 & 1 & 0 & 0 & 0\\
1 & 0 & 1 & 1 & 0 & 0 & 0 & 0\\
1 & 1 & 0 & 1 & 1 & 0 & 0 & 0\end{array}\right)$.
\end{center}

Similarly, for cylindrical quadratic lattice graph $H_{n,m}$, its
transfer matrix $T_{n}$ is an
$(F_{n-1}+F_{n+1})\times(F_{n-1}+F_{n+1})$ symmetric matrix of $0$'s
and $1$'s, defined as follows. The rows and columns of $T_{m}$ are
indexed by vectors of $C_{n}$, and the entry of $T_{n}$ in position
$(\alpha,\beta)$ is $1$ if the vectors $\alpha,\beta$ are
orthogonal, and is $0$ otherwise. For $n=4$ the possible row vectors
of $C{}_{n}$ are $(0,0,0,0)$, $(1,0,0,0)$, $(0,1,0,0)$, $(0,0,1,0)$,
$(0,0,0,1)$, $(1,0,1,0)$, $(0,1,0,1)$. If we index the row and
column in this order, then the transfer matrix of $H_{n,m}$ is

\begin{center}
$B_{4}$=$\left(\begin{array}{ccccccc}
1 & 1 & 1 & 1 & 1 & 1 & 1\\
1 & 0 & 1 & 1 & 0 & 1 & 1\\
1 & 1 & 0 & 1 & 1 & 1 & 0\\
1 & 1 & 1 & 0 & 0 & 1 & 1\\
1 & 0 & 1 & 0 & 0 & 1 & 1\\
1 & 1 & 1 & 1 & 1 & 0 & 0\\
1 & 1 & 0 & 1 & 1 & 0 & 0\end{array}\right).$
\end{center}

 For a plane
quadratic lattice graph $G_{m,n}$, by adding two crossed diagonals
in each square inner face we get a crossed quadratic lattice graph,
which will be denoted as $\check G_{m,n }$. The transfer matrix of
$\check G_{m,n }$ can be defined in a similar way.  We say that
$G_{m,n}$ and $\check G_{m,n }$ are lattice graphs  with transfer
multiplicity one, since  computing the number of independent sets
for each of them only needs to employ one transfer matrix. This may
not hold for other lattices. For example, in \cite{zz}, for an
aztec diamond we need to introduce two transfer matrices, each of
which is the transpose of the other. So the transfer multiplicity of
an  aztec diamond is two. There are also lattices with transfer
multiplicity three; one example is the 8.8.4 lattice that will be
discussed later in section 4.

The entropy constant of a plane lattice is defined by
$\eta=\underset{m,n\rightarrow\infty}{\lim}f(m,n)^{1/k(m,n)}$ where
$f(m,n)$ denotes the number of independent sets of $G_{m,n}$ and
$k(m,n)$ denotes the number of vertices of $G_{m,n}$. The entropy
constants of cylindrical and toridal lattices can be defined
similarly. As in \cite{cw}, we have $f(m,n)=\underset{u,v\in
P_{m}}{\sum}T_{u,v}^{n}=1\cdot T^{n}1$ for $G_{m,n}$ and $H_{n,m}$.
Clearly, Trace$(T^{n})$=$\underset{u\in P_{m}}{\sum} T_{u,u}^{n}.$

 In \cite{cw}, Calkin and Wilf proved the existence of the entropy
constant of plane quadratic lattice and established its upper and
lower bounds. Two natural problems are to consider the entropy
constants for lattices on cylinder or torus. Note that the method
of Calkin and Wilf's is valid for the lattices with the same
symmetric transfer matrices in both horizontal and vertical
directions. In this paper, we will consider three types of $2$-dimensional
lattices: plane quadratic lattice, generalized aztec diamonds lattice and
8.8.4 lattice. We will also consider crossed quadratic lattice which is
a non-planar lattice. We will show that for each type of these four lattices, the
entropy constant is the same no matter the lattice is on plane,
cylinder or torus. Furthermore, the upper and lower bounds of the entropy constant will be established for
crossed quadratic lattice, generalized aztec diamonds lattice and 8.8.4 lattice.

\section{Lattices with transfer multiplicity one}

\begin{figure}[htbp]
  \centering
 \scalebox{0.5}{\includegraphics{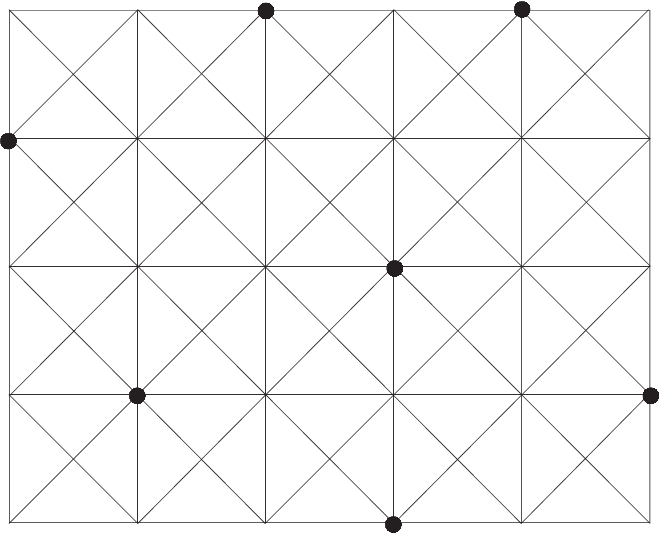}}
  \caption{}\label{fig:fig2.1}
\end{figure}

 The entropy constant of the plane quadratic lattice was already
discussed in \cite{cw}. Now we consider the crossed quadratic
lattice, which is
 obtained from the plane quadratic
lattice  by  adding two crossed diagonals to each square inner face.
  Fig.\ref{fig:fig2.1} shows the  crossed quadratic lattice
 graph $\check G_{4,5}$ that is a  part  of the crossed quadratic lattice,
 where an independent set $S$ is indicated by small circles.
The portions of S that lie in each of the columns can be represented
by $5$-vectors $(0,1,0,0,0)$, $(0,0,0,0,1)$, $(1,0,1,0,0)$,
$(0,0,0,0,1)$, $(1,0,1,0,0)$, $(0,0,0,0,1)$.

%%%%%%%%%%%%%%%%%%%%%%%%%%%%%%%%%%%%%%%%%%
%%%%%%%%%%%%%
%%%%%%%%%%%%%%%%%%%%%%%%%%%%%%%%%%%%%%%%%%

%\begin{figure}[h]
%\center
%  \includegraphics[scale=0.5]{./eps/fig2.1.eps}\\
%  \caption{}\label{fig:fig2.1}
%\end{figure}

Note that by allowing the crossed edges, the plane crossed  lattice are no longer planar
graphs. However,
 it is easy to see
that the transfer matrix of the plane crossed lattice graph $\check
G_{m,n }$ is an $F_{m+3}\times F_{m+3}$ matrix of 0's and 1's. When
$m=3$, the possible column vectors of $P_{m}$ are $(0,0,0,0)$,
$(1,0,0,0)$, $(0,1,0,0)$, $(0,0,1,0)$, $(0,0,0,1)$, $(1,0,1,0)$,
$(1,0,0,1)$, $(0,1,0,1)$. If we index the rows and columns in this
order, then the transfer matrix of $\check G_{m,n }$ is

\begin{center}
$T_{3}=\left(\begin{array}{cccclccc}
1 & 1 & 1 & 1 & 1 & 1 & 1 & 1\\
1 & 0 & 0 & 1 & 1 & 0 & 0 & 0\\
1 & 0 & 0 & 0 & 1 & 0 & 0 & 0\\
1 & 1 & 0 & 0 & 0 & 0 & 0 & 0\\
1 & 1 & 1 & 0 & 0 & 0 & 0 & 0\\
1 & 0 & 0 & 0 & 0 & 0 & 0 & 0\\
1 & 0 & 0 & 0 & 0 & 0 & 0 & 0\\
1 & 0 & 0 & 0 & 0 & 0 & 0 & 0\end{array}\right)$.
\end{center}
   Clearly,  the plane crossed lattice is a non-planar lattice
 with transfer multiplicity one.

For a given $\check G_{m,n }$, by identifying its left edges with
the right edges correspondingly, we get the cylindrical crossed
lattice graph $\check H_{n,m}$.
 The transfer matrix of
 $\check H_{n,m}$ is an $(F_{n-1}+F_{n+1})\times(F_{n-1}+F_{n+1})$ symmetric matrix of $0$'s
and $1$'s. When $n=4$ the possible row vectors of its $C_{n}$ are
$(0,0,0,0)$, $(1,0,0,0)$, $(0,1,0,0)$, $(0,0,1,0)$, $(0,0,0,1)$,
$(1,0,1,0)$, $(0,1,0,1)$. If we index the rows and columns in this
order, then the transfer matrix of $\check H_{n,m}$ is

\begin{center}
$B_{4}=\left(\begin{array}{ccccccc}
1 & 1 & 1 & 1 & 1 & 1 & 1\\
1 & 0 & 0 & 1 & 0 & 0 & 0\\
1 & 0 & 0 & 0 & 1 & 0 & 0\\
1 & 1 & 0 & 0 & 0 & 0 & 0\\
1 & 0 & 1 & 0 & 0 & 0 & 0\\
1 & 0 & 0 & 0 & 0 & 0 & 0\\
1 & 0 & 0 & 0 & 0 & 0 & 0\end{array}\right)$.
\end{center}

In the following lemma, we will establish a relation between the
transfer matrix of $G_{m,n}$ (cylindrical lattice $H_{n,m}$) and the
number of independent sets of the corresponding cylinder lattice
graph $H_{n,m}$ (toroidal lattice graph $S_{n,m}$). Then we will
prove our main theorem and provide numerical upper and lower bounds
of the entropy constants for the quadratic and crossed lattices on
cylinder and torus.

\begin{lemma}
For each plane lattice graph $G_{m,n}$ with transfer matrix $T_m$,
transfer multiplicity one and a positive integer $n$, the trace of
$T_{m}^{n}$ is equal to the number of independent sets of the
corresponding cylinder lattice graph $H_{n,m}$. For each cylinder
lattice graph $H_{n,m}$ with transfer matrix $T_n$, transfer
multiplicity one and a positive integer $m$, the trace of
$T_{n}^{m}$ is equal to the number of independent sets of the
corresponding torus lattice graph $S_{n,m}$.

\end{lemma}
\begin{proof}
Recall that $H_{n,m}$ can be obtained by identifying the left column
and the right column of $G_{m,n}$. Thus there is a bijection between
the independent sets of $H_{n,m}$ and the independent sets of
$G_{m,n}$ whose left and right column vector are the same. And the
latter is the trace of $T_{m}^{n}$.

Similarly, $S_{n,m}$ can be obtained by identifying the top row and
the bottom row of $H_{n,m}$. Thus there is a bijection between the
independent sets of $S_{n,m}$ and the  independent sets of $H_{n,m}$
whose corresponding top and bottom row vectors are the same. And the
latter is the trace of $T_{n}^{m}$.
\end{proof}

\begin{theorem}
For a lattice with transfer multiplicity one, if in both directions
the transfer matrices are the same real symmetric matrix, then its
entropy constants on plane, cylinder and torus are the same.
\end{theorem}
\begin{proof}
Let $T$ be the transfer matrix of columns of $G_{m,n}$ with  the
characteristic polynomial $f(x)$,  and let
\begin{center}
$\mathbf{1}$=$(1,1,\cdots,1)_{1\times g(m)}$
\end{center}
where $g(m)$ equals  the number of independent sets on the first
column of $G_{m,n}$. Particularly, for the quadratic lattice and the
crossed lattice, $g(m)=F_{m+3}.$

By Hamilton-Cayley Theorem we have\\
\begin{center}
$f(T)=a_{0}I+a_{1}T+a_{2}T^{2}+\cdots+a_{g(m)}T^{g(m)}=0$.
\end{center}

Put $f_{0}=\mathbf{1}$ and $f_{n+1}=Tf_{n}$, then
\begin{center}
$f_{n}=T^{n}\mathbf{1}$.
\end{center}

Let
\begin{center}
$u_{m,n}=\mathbf{\mathbf{1}}\cdot f_{n}=\mathbf{1}\cdot
T^{n}\mathbf{1}$,
\end{center}
where $u_{m,n}$ is the sum of all entries of the matrix $T^{n}$,
namely the number of independent sets of $G_{m,n}$.

Thus \\
\begin{center}
$\mathbf{1}\cdot
f(T)\mathbf{1}=a_{0}u_{m,0}+a_{1}u_{m,1}+a_{2}u_{m,2}+\cdots+a_{g(m)}u_{m,g(m)}=0$.
\end{center}

In general $u_{m,n}$ satisfies the following recurrence relation:\\
\begin{center}
$a_{0}u_{m,n}+a_{1}u_{m,n+1}+a_{2}u_{m,n+2}+\cdots+a_{g(m)}u_{m,n+g(m)}=0$.
\end{center}

By a well known theorem on difference equations (see, for example,
\cite{h}), the
characteristic polynomial of this linear recurrence relation can be written as\\

\begin{center}
$f(x)=(x-\lambda_{1})^{e_{1}}(x-\lambda_{2})^{e_{2}}\cdots(x-\lambda_{s})^{e_{s}}$
where $\lambda_{1}\leq\lambda_{2}\leq\cdots\leq\lambda_{s}$.
\end{center}

Note that $\lambda_{s}$ depends on m. Then\\
\begin{center}
$u_{m,n}=\overset{}{\underset{}{\overset{s}{\underset{i=1}{\sum}}p_{i}(n)\lambda_{i}^{n}}}$
\end{center}

\noindent for all $n$ where $p_{i}(n)$ is a polynomial  with degree
at most $e_{i}-1$ in  $n$. The coefficients $p_{i}(n)$ of the
polynomial $u_{m,n}$ are determined by the initial values
\begin{center}
$u_{m,0},u_{m,1},\cdots,u_{m,g(m)-1}.$
\end{center}

Note that  the first row and the first column of $T$ are both
vectors of all $1$'s. Thus the matrix $T$ is non-negative,
irreducible (its corresponding digraph is strongly connected) and
prime. Hence the spectral radius of $T$ is a simple positive real
eigenvalue $\lambda_{s}$ with magnitude greater than any other
eigenvalues'(see \cite{v})
and $p_{s}(n)$ is a positive constant. It is not difficult to see that\\
\begin{center}
$\underset{n\rightarrow\infty}{\lim}(u_{m,n})^{1/n}=\lambda_{s}.$
\end{center}

Since in both directions the real and symmetric transfer matrices
are the same, by the same method in \cite{cw}, one can show that the
following double
limit exists:\\

\begin{center}
$\eta_{1}=\underset{m,n\rightarrow\infty}{\lim}(u_{m,n})^{1/m}=\underset{m\rightarrow\infty}{\lim}\lambda_{s}^{m}$.
\end{center}

Now by Lemma 2.1 we can see that the trace of the $n$-th power of
the transfer matrix of $G_{m,n}$  is equal to $v_{m,n}$, the number
of independent sets of $H_{n,m}$.

Since the trace of $T^{n}$ is the sum of its eigenvalue of $T^{n}$
and each one is an $n$-th power of an eigenvalue of $T$, we see that

if
\begin{center}
$v_{m,1}=Trace(T)=\overset{}{\overset{s}{\underset{i=1}{\sum}}e_{i}\lambda_{i},}$
 \end{center}
then
\begin{center}
$v_{m,n}=Trace(T^{n})=\overset{}{\underset{i=1}{\overset{s}{\sum}}e_{i}\lambda_{i}^{n}}$,
\end{center}

where $e_{i}$ is the multiplicity of eigenvalue $\lambda_{i}$.

Since $\lambda_{s}$ is the simple largest eigenvalue, then

\begin{center}
$\underset{n\rightarrow\infty}{\lim}(v_{m,n})^{1/n}=\lambda_{s}$.
\end{center}

So, the entropy constant of $H_{n,m}$

\begin{center}
$\eta_{2}=\underset{m,n\rightarrow\infty}{\lim}(v_{m,n})^{1/m(n-1)}
=\underset{m\rightarrow\infty}{\lim}\lambda_{s}^{1/m}=\eta_{1}$,
\end{center}

That is, the entropy constant of $G_{m,n}$ and $H_{n,m}$ are the
same. Similarly, using the second conclusion of Lemma 2.1 we can
prove that entropy constants of $H_{n,m}$ and $S_{n,m}$ are the
same. This completes the proof.
\end{proof}

In \cite{cw} Calkin and Wilf already obtained a good estimate for
the entropy constant of the plane quadratic lattice. Now, by Theorem
2.2 we immediately see that the entropy constants of the quadratic
lattice on plane, cylinder and torus are all between
$1.503047782...$ and $1.5035148...$.

Since the transfer matrix of crossed quadratic lattice is symmetric,
the same approach in \cite{cw} and the proof of theorem 2.1 can be
processed here. Thus we can use the same method in \cite{cw} to get
the upper and lower bound of the entropy constant of crossed
quadratic lattice. The lower bound of
$\underset{m,n\rightarrow\infty}{\lim}f(m,n)^{1/mn}$ is
$(\frac{\lambda_{p+2q}}{\lambda_{2q}})^{1/p}$ where $\lambda$'s are
the largest eigenvalues of corresponding $T$'s and the upper bound
of $\underset{m,n\rightarrow\infty}{\lim}f(m,n)^{1/mn}$ is
$(\xi_{2k})^{1/2k}$ where $\xi$'s are the largest eigenvalues of
corresponding $B$'s. Let $p=4$, $q=4$ and $k=6$, we have
$1.342542258...\leq\eta\leq1.342652572...$.
\section{Lattices with transfer multiplicity two}

In this section we consider the lattices with transfer multiplicity
two, i.e., the lattices for which we need two transfer matrices to
compute the  number of independent sets of the lattice graphs. One
of such lattices is
 inspired by the famous aztec diamonds. The study of enumeration of perfect matchings,
spanning trees and independent sets of an aztec diamond can be found
in \cite{c,mc,k} and the references cited therein.

\begin{figure}[htbp]
  \centering
 \scalebox{0.5}{\includegraphics{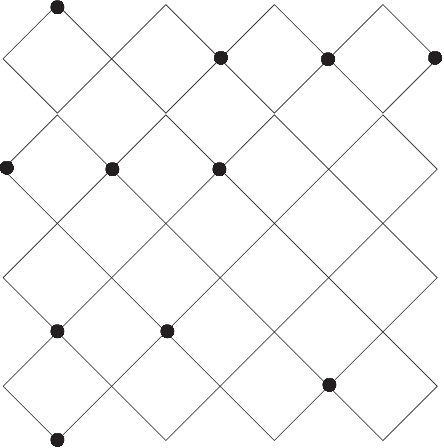}}
  \caption{}\label{fig:fig3.1}
\end{figure}

Let $L_{i}$ be the path with $i$ vertices $1,2$,$\cdots$,$i$. The
tensor product of two paths $L_{m}\otimes L_{n}$ is the graph on
$m\times n$ vertices $\left\{ \left(x,y\right):1\leq x\leq m,1\leq
y\leq n\right\} $, with $\left(x,y\right)$ adjacent to
$(x^{'},y^{'})$ if and only if $|x-x^{'}|=|y-y^{'}|=1$. This graph
consists of two connected components, the one with the vertices
\{$(x,y)|x+y$ is odd\}, denoted $O(L_{2m+1}\otimes L_{2m+1})$, is
called the aztec diamond of order m. More generally
$O(L_{2m+1}\otimes L_{2n+1})$ is called the generalized aztec
diamond of order $m\times n$ introduced by the author in \cite{zz}.
An independent set of a generalized aztec diamond can be represented
by an ordered list of column vectors. Fig.\ref{fig:fig3.1} shows an
example. For $O(L_{9}\otimes L_{9})$, the given $S$ (the big dots)
can be represented by the ordered list of $9$ column vectors:
$(0,1,0,0)$, $(1,0,0,1,1)$, $(0,1,0,0)$, $(0,0,0,1,0)$, $(1,1,0,0)$,
$(0,0,0,0,0)$, $(0,0,0,1)$, $(0,1,1,0,0)$, $(0,0,0,1)$,
$(1,1,1,0,0)$.

%%%%%%%%%%%%%%%%%%%%%%%%%%%%%%%%%%%%%%%%%%
%%%%%%%%%%%%%
%%%%%%%%%%%%%%%%%%%%%%%%%%%%%%%%%%%%%%%%%%
%\begin{figure}[h]
%\center
%  \includegraphics[scale=0.5]{./eps/fig3.1.eps}\\
%  \caption{}\label{fig:fig3.1}
%\end{figure}
%\includegraphics{\string"aztec diamond\string".png}
%\noindent
By identifying the top row and the bottom row of a generalized aztec
diamond, we get  a cylindrical generalized aztec diamond. And, the
toroidal generalized aztec diamonds
  can be obtained by
identifying the left cycle and the right cycle of a cylindrical
generalized aztec diamond.

Now we consider the transfer matrix of $O(L_{2m+1}\otimes
L_{2n+1})$. It is clear that the generation of each independent set
of $O(L_{2m+1}\otimes L_{2n+1})$ involves $2n$ assembling steps
after the first column $v_{1}$ is established (For simplicity, the
assembling of the $(i+1)$-th column $v_{i+1}$ to the $ith$ column
$v_{i}$ is called step $i$.)  Step one is to assemble $v_{2}$ to the
right side of $v_{1}$. The transfer matrix representing step one,
denoted $T_{m_{1}}$, can be constructed as follows. Let $R_{m}$ be
the set of all possible vectors $v_{1}$. Obviously, $R_{m}$ consists
of all m-dimensional vectors of $0$'s and $1$'s, so it has $2^{m}$
vectors. Similarly, the set of all possible $v_{2}$ is the set
$R_{m+1}$ of all $(m+1)$-dimensional vectors of $0$'s and $1$'s, and
$R_{m+1}$ has $2^{m+1}$ vectors. Then the transfer matrix
$T_{m_{1}}=(T_{v_{1},v_{2}})$ is a $2^{m}\times2^{m+1}$ matrix whose
rows are indexed by vectors of $R_{m}$ and columns are indexed by
vectors of $R_{m+1}$, where $T_{v_{1},v_{2}}=1$ if $v_{1}$ and
$v_{2}$ represent possible consecutive pair of columns in an
independent set of $O(L_{2m+1}\otimes L_{2n+1})$, and
$T_{v_{1},v_{2}}=0$ otherwise. Similarly, the transfer matrix for
step two is a $2^{m+1}\times2^{m}$ matrix $T_{m_{2}}$ which is the
transpose of $T_{m_{1}}$. It is easily seen that $T_{m_{1}}$ is the
transfer matrix for every step $i$ where $i$ is odd and $T_{m_{2}}$
is the transfer matrix for every step $i$ where $i$ is even. Thus,
if we take the transfer matrix of generalized aztec diamond to be
$T_{m}=T_{m_{1}}T_{m_{2}}$ then it is a $2^{m}\times2^{m}$ symmetric
matrix. Furthermore, the transfer multiplicity of $T_{m}$ can be
considered as one and the results of Theorem 2.2 hold.

When $m=3$, if we index the rows and columns in increasing order in
binary numbers, then the transfer matrix of $G_{m,n}$ is an
$8\times8$ matrix given by the product of $T_{3_{1}}$ and
$T_{3_{2}}$. Here
\begin{center}
$T_{3_{1}}=\left(\begin{array}{cccccccccccccccc}
1 & 1 & 1 & 1 & 1 & 1 & 1 & 1 & 1 & 1 & 1 & 1 & 1 & 1 & 1 & 1\\
1 & 0 & 0 & 0 & 1 & 0 & 0 & 0 & 1 & 0 & 0 & 0 & 1 & 0 & 0 & 0\\
1 & 1 & 0 & 0 & 0 & 0 & 0 & 0 & 1 & 1 & 0 & 0 & 0 & 0 & 0 & 0\\
1 & 0 & 0 & 0 & 0 & 0 & 0 & 0 & 1 & 0 & 0 & 0 & 0 & 0 & 0 & 0\\
1 & 1 & 1 & 1 & 0 & 0 & 0 & 0 & 0 & 0 & 0 & 0 & 0 & 0 & 0 & 0\\
1 & 0 & 0 & 0 & 0 & 0 & 0 & 0 & 0 & 0 & 0 & 0 & 0 & 0 & 0 & 0\\
1 & 1 & 0 & 0 & 0 & 0 & 0 & 0 & 0 & 0 & 0 & 0 & 0 & 0 & 0 & 0\\
1 & 0 & 0 & 0 & 0 & 0 & 0 & 0 & 0 & 0 & 0 & 0 & 0 & 0 & 0 &
0\end{array}\right)$,
\end{center}
\begin{center}
$T_{3_{2}}=\left(\begin{array}{cccccccc}
1 & 1 & 1 & 1 & 1 & 1 & 1 & 1\\
1 & 0 & 1 & 0 & 1 & 0 & 1 & 0\\
1 & 0 & 0 & 0 & 1 & 0 & 0 & 0\\
1 & 0 & 0 & 0 & 1 & 0 & 0 & 0\\
1 & 1 & 0 & 0 & 0 & 0 & 0 & 0\\
1 & 0 & 0 & 0 & 0 & 0 & 0 & 0\\
1 & 0 & 0 & 0 & 0 & 0 & 0 & 0\\
1 & 0 & 0 & 0 & 0 & 0 & 0 & 0\\
1 & 1 & 1 & 1 & 0 & 0 & 0 & 0\\
1 & 0 & 1 & 0 & 0 & 0 & 0 & 0\\
1 & 0 & 0 & 0 & 0 & 0 & 0 & 0\\
1 & 0 & 0 & 0 & 0 & 0 & 0 & 0\\
1 & 1 & 0 & 0 & 0 & 0 & 0 & 0\\
1 & 0 & 0 & 0 & 0 & 0 & 0 & 0\\
1 & 0 & 0 & 0 & 0 & 0 & 0 & 0\\
1 & 0 & 0 & 0 & 0 & 0 & 0 & 0\end{array}\right)$.
\end{center}

Thus,
\begin{center}
$T_{3}=$$\left(\begin{array}{cccccccc}
16 & 4 & 4 & 2 & 4 & 1 & 2 & 1\\
4 & 4 & 2 & 2 & 1 & 1 & 1 & 1\\
4 & 2 & 4 & 2 & 2 & 1 & 2 & 1\\
2 & 2 & 2 & 2 & 1 & 1 & 1 & 1\\
4 & 1 & 2 & 1 & 4 & 1 & 2 & 1\\
1 & 1 & 1 & 1 & 1 & 1 & 1 & 1\\
2 & 1 & 2 & 1 & 2 & 1 & 2 & 1\\
1 & 1 & 1 & 1 & 1 & 1 & 1 & 1\end{array}\right)$.
\end{center}

Let $R_{n}$ denote all $n$-dimensional vectors of $0$'s and $1$'s.
The cylindrical generalized aztec diamond, obtained by identifying
the left and right columns of the generalized aztec diamond
$O(L_{2m+1}\otimes L_{2n+1})$, can also be seen as obtained by
beginning with some vector in $R_{n}$, gluing to the top  a new one
in $R_{n}$ such that the vertices represented by the 1's in these
two vectors are not adjacent until $2m+1$ vectors are glued. It is
clear that for each $i$, step $i$ and step $i+2$ can be represented
by the same transfer matrix. Let  $B_{n_{1}}$ ($B_{n_{2}}$, resp)
denote the transfer matrix  of the cylindrical generalized aztec
diamond,
 which represents every step $i$ for $i$ odd (even, resp).  Clearly,
  $B_{n_{1}}$ is a
$2^{n}\times2^{n}$ matrix of 0's and $1$'s whose rows and columns
are indexed by vectors of $R_{n}$, the entry of $B_{n_{1}}$ in
position $(\alpha,\beta)$ is $1$ if the vectors represent possible
consecutive pair of rows in an independent set of $O(L_{2m+1}\otimes
L_{2n+1})$ on cylinder, and is $0$ otherwise. It is no difficult to
see that $B_{n_{2}}$ is the transpose of $B_{n_{1}}$. Thus  we get
the transfer matrix of cylindrical aztec diamond as
$B_{n}=B_{n_{1}}B_{n_{2}}$, which is a $2^{n}\times2^{n}$ symmetric
matrix. Furthermore, the transfer multiplicity of $B_{n}$ is one.

Take $n=3$ as an example. If we index the rows and columns in
increasing order in binary numbers, then the transfer matrix of
$H_{n,m}$ is an $8\times8$ matrix given by the product of
$B_{3_{1}}$ and $B_{3_{2}}$, \noindent where
\begin{center}
$B_{3_{1}}=\left(\begin{array}{cccccccc}
1 & 1 & 1 & 1 & 1 & 1 & 1 & 1\\
1 & 0 & 1 & 0 & 0 & 0 & 0 & 0\\
1 & 0 & 0 & 0 & 1 & 0 & 0 & 0\\
1 & 0 & 0 & 0 & 0 & 0 & 0 & 0\\
1 & 1 & 0 & 0 & 0 & 0 & 0 & 0\\
1 & 0 & 0 & 0 & 0 & 0 & 0 & 0\\
1 & 0 & 0 & 0 & 0 & 0 & 0 & 0\\
1 & 0 & 0 & 0 & 0 & 0 & 0 & 0\end{array}\right)$
\end{center}
 and
\begin{center}
$B_{3_{2}}=\left(\begin{array}{cccccccc}
1 & 1 & 1 & 1 & 1 & 1 & 1 & 1\\
1 & 0 & 0 & 0 & 1 & 0 & 0 & 0\\
1 & 1 & 0 & 0 & 0 & 0 & 0 & 0\\
1 & 0 & 0 & 0 & 0 & 0 & 0 & 0\\
1 & 0 & 1 & 0 & 0 & 0 & 0 & 0\\
1 & 0 & 0 & 0 & 0 & 0 & 0 & 0\\
1 & 0 & 0 & 0 & 0 & 0 & 0 & 0\\
1 & 0 & 0 & 0 & 0 & 0 & 0 & 0\end{array}\right)$.
\end{center}

\noindent Thus

\begin{center}
$B_{3}=\left(\begin{array}{cccccccc}
8 & 2 & 2 & 1 & 2 & 1 & 1 & 1\\
2 & 2 & 1 & 1 & 1 & 1 & 1 & 1\\
2 & 1 & 2 & 1 & 1 & 1 & 1 & 1\\
1 & 1 & 1 & 1 & 1 & 1 & 1 & 1\\
2 & 1 & 1 & 1 & 2 & 1 & 1 & 1\\
1 & 1 & 1 & 1 & 1 & 1 & 1 & 1\\
1 & 1 & 1 & 1 & 1 & 1 & 1 & 1\\
1 & 1 & 1 & 1 & 1 & 1 & 1 & 1\end{array}\right)$.
\end{center}

Since the generalized aztec diamond lattice has the same symmetric
transfer matrices in both horizontal and vertical directions, the
same approach in \cite{cw} can be taken here. So, the lower bound of
$\underset{m,n\rightarrow\infty}{\lim}f(m,n)^{1/mn}$ is
$(\frac{\lambda_{p+2q}}{\lambda_{2q}})^{1/p}$ where $\lambda$'s are
the largest eigenvalues of corresponding $T$'s, and the upper bound
of $\underset{m,n\rightarrow\infty}{\lim}f(m,n)^{1/mn}$ is
$(\xi_{2k})^{1/2k}$ where $\xi$'s are the largest eigenvalues of
corresponding $B$'s. Taking $p=2$ , $q=4$ and $k=5$, we get
\begin{center}
$2.259132578...\leq\underset{m,n\rightarrow\infty}{\lim}f(m,n)^{1/mn}\leq2.259154406...$.
\end{center}

Note that the entropy constant of the generalized aztec diamond
lattice is

$\eta=\underset{m,n\rightarrow\infty}{\lim}f(m,n)^{1/(2mn+m+n)}=\underset{m,n\rightarrow\infty}{\lim}f(m,n)^{1/2mn}$
where $2mn+m+n$ is the number of vertices of vertices of generalized
aztec diamond. Then we see that  the entropy constant of the
generalized aztec diamond lattice
 is
$(\underset{m,n\rightarrow\infty}{\lim}f(m,n)^{1/mn}){}^{1/2}$,
which is between $1.503041110...$ and $1.503048371...$.

%%%%%%%%%%%%%%%%%%%%%%%%%%%%%%%%%%%%%%%%%%%%%%%%%%%%%%%%%%%%%%%%%%%%%%%%

\section{Lattices with transfer multiplicity three}
\begin{figure}[htbp]
  \centering
 \scalebox{0.7}{\includegraphics{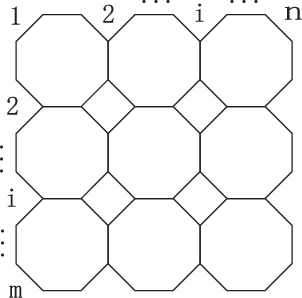}}
  \caption{}\label{fig:fig4.1}
\end{figure}

The 8.8.4 lattice graphs $G_{m,n}$, as shown in
Fig.\ref{fig:fig4.1}, are finite subgraphs of the 8.8.4 tiling of
Euclidean plane. Some study of the properties of 8.8.4 lattice
graphs can be found in \cite{sn,yyz}.

%%%%%%%%%%%%%%%%%%%%%%%%%%%%%%%%%%%%%%%%%%
%%%%%%%%%%%%%
%%%%%%%%%%%%%%%%%%%%%%%%%%%%%%%%%%%%%%%%%%

%\begin{figure}[h]
%\center
%  \includegraphics[scale=0.8]{./eps/fig4.1.eps}\\
%  \caption{}\label{fig:fig4.1}
%\end{figure}

The cylindrical 8.8.4 lattice graphs  are obtained by identifying
the top row and the bottom row of the 8.8.4 lattice graphs
$G_{m,n}$. The toroidal 8.8.4 lattice graphs  can be obtained by
identifying the left  cycle and the right cycle of the cylindrical
8.8.4 lattice graphs.

Consider the transfer matrix of the 8.8.4 lattice graph $G_{m,n}$.
Define the assembling of the $(i+1)th$ column $v_{i+1}$ to the $ith$
column $v_{i}$ as step $i$ as we did before. We can see that the
generation of each independent set of $G_{m,n}$ involves $3n-3$
steps after the first column $v_{1}$ determined. Step one is to
assemble $v_{2}$ to the right side of $v_{1}$. The transfer matrix
representing step one, denoted $T_{m_{1}}$, can be constructed as
follows.
  Let $R_{m}$ denote the set of all
possible vectors  which can appear as $v_{1}$.  Clearly $R_{m}$
consists of $(2m-2)$-vectors of $0$'s and $1$'s in which no
consecutive 1's occupy  the positions of the $2k-1$-th and $2k$-th
entries, for  $1\leq k\leq m-1$. Since there are three possibilities
at each pair of consecutive $2k-1$-th and $2k$-th $1$ positions, the
set $R_{m}$ has $3^{m-1}$ vectors. The set of all possible vectors
$v_{2}$, denoted as $K_{m}$, consists of all $m$-vectors of $0$'s
and $1$'s. So $K_{m}$ has $2^{m}$ vectors. Thus the transfer matrix
$T_{m_{1}}=(T_{v_{1},v_{2}})$ is a $3^{m-1}\times2^{m}$ matrix whose
rows are indexed by vectors of $R_{m}$ and columns are indexed by
vectors of $K_{m}$, where $T_{v_{1},v_{2}}=1$ if $v_{1}$ and $v_{2}$
represent possible consecutive pair of columns in an independent set
of $G_{m,n}$ and $T_{v_{1},v_{2}}=0$ otherwise. Similarly, the
transfer matrix for step three is a $2^{m}\times3^{m-1}$ matrix
$T_{m_{3}}$ which is the transpose of $T_{m_{1}}$, and the transfer
matrix for step two is a $2^{m}\times2^{m}$ matrix $T_{m_{2}}$ whose
rows and columns are indexed by vectors in $K_{m}$, $T_{m_{2}}$'s
entry in position $(\alpha,\beta)$ is $1$ if the vectors represented
by $\alpha$,$\beta$ are  orthogonal, and is $0$ otherwise. Note that
$T_{m_{1}}$ is the transfer matrix for every step $i$ when
$i=3k+1$$(0\leq k\leq n-2)$, $T_{m_{2}}$ is the transfer matrix for
every step $i$ when $i=3k+2$$(0\leq k\leq n-2)$ and $T_{m_{3}}$ is
the transfer matrix for every step $i$ when $i=3k$ $(0\leq k\leq
n-1)$. Thus, the transfer multiplicity of $G_{m,n}$ is 3. Thus, if
we take the transfer matrix of $G_{m,n}$ to be
$T_{3}$=$T_{3_{1}}$$T_{3_{2}}T_{3_{3}}$, then it is a
$2^{m}\times2^{m}$ symmetric matrix with transfer multiplicity one.

When $m=2$, if we index the rows and columns in increasing order in
binary numbers, then the transfer matrix of $G_{m,n}$ is a
$3\times3$ matrix given by the product of $T_{2_{1}}$, $T_{2_{2}}$
and $T_{2_{3}}$.

\begin{center}
$T_{2_{1}}=\left(
              \begin{array}{cccc}
                1 & 1 & 1 & 1 \\
                1 & 0 & 1 & 0 \\
                1 & 1 & 0 & 0 \\
              \end{array}\right)$
\end{center}

\begin{center}
$T_{2_{2}}=\left(
              \begin{array}{cccc}
                1 & 1 & 1 & 1 \\
                1 & 0 & 1 & 0 \\
                1 & 1 & 0 & 0 \\
                1 & 0 & 0 & 0 \\
              \end{array}\right)$.
\end{center}

\begin{center}
$T_{2_{3}}=\left(\begin{array}{ccc}
                                 1 & 1 & 1 \\
                                 1 & 0 & 1 \\
                                 1 & 1 & 0 \\
                                 1 & 0 & 0\end{array}\right)$
\end{center}

\begin{center}
Thus $T_{2}$=$\left(
                            \begin{array}{ccc}
                              9 & 6 & 6 \\
                              6 & 3 & 4 \\
                              6 & 4 & 3 \\\end{array}\right)$.
\end{center}

Consider the transfer matrix of $H_{n,m}$ which is obtained from
$G_{m,n}$ by identifying its left column and  right column. The
transfer matrix $B_{n_{1}}$, which represents every $(3k+1)$-th
$(0\leq k\leq n-2)$ step, can be defined as a $3^{n-1}\times2^{n-1}$
matrix of $0$'s and $1$'s as follows. The rows of $B_{n_{1}}$ are
indexed by vectors of $R_{n}$ and columns are indexed by vectors of
$K_{n-1}$, and the entry of $B_{n_{1}}$ in position $(\alpha,\beta)$
is $1$ if  $\alpha,\beta$ represent
 possible consecutive pair of rows
in an independent set of $H_{n,m}$, and is $0$ otherwise. Let
$B_{n_{3}}$ denote the transfer matrix that represents every $3k$-th
$(0\leq k\leq n-1)$ step.  It is no difficult to see that
$B_{n_{3}}$ is the transpose of $B_{n_{1}}$. The transfer matrix
$B_{n_{2}}$ that represents every $(3k+2)$-th $(1\leq k\leq n-2)$
step is a $2^{n-1}\times2^{n-1}$ matrix whose rows and columns are
indexed by vectors of $K_{n-1}$. The entry of $B_{n_{2}}$ in
position $(\alpha$,$\beta)$ is $1$ if $\alpha$,$\beta$ are
orthogonal, and is $0$ otherwise. Thus if we take the transfer
matrix of $H_{n,m}$ to be $B_{n}=B_{n_{1}}B_{n_{2}}B_{n_{3}}$, then
it is a $2^{n-1}\times2^{n-1}$ symmetric matrix with transfer
multiplicity one.

When $n=3$, if we index the rows and columns in increasing order in
binary numbers, then the transfer matrix of $H_{n,m}$ is a
$8\times8$ matrix given by product of $B_{3_{1}}$, $B_{3_{2}}$ and
$B_{3_{3}}$, where
\begin{center}
$B_{3_{1}}=\left(
              \begin{array}{cccc}
                1 & 1 & 1 & 1 \\
                1 & 1 & 0 & 0 \\
                1 & 0 & 1 & 0 \\
                1 & 0 & 1 & 0 \\
                1 & 0 & 0 & 0 \\
                1 & 0 & 1 & 0 \\
                1 & 1 & 0 & 0 \\
                1 & 1 & 0 & 0 \\
                1 & 0 & 0 & 0 \\
              \end{array}\right)$,
\end{center}

\begin{center}
$B_{3_{2}}=\left(\begin{array}{cccc}
1 & 1 & 1 & 1\\
1 & 0 & 1 & 0\\
1 & 1 & 0 & 0\\
1 & 0 & 0 & 0\end{array}\right)$.
\end{center}

Thus
\begin{center}
$B_{3_{3}}=\left(
             \begin{array}{ccccccccc}
               1 & 1 & 1 & 1 & 1 & 1 & 1 & 1 & 1 \\
               1 & 1 & 0 & 0 & 0 & 0 & 1 & 1 & 0 \\
               1 & 0 & 1 & 1 & 0 & 1 & 0 & 0 & 0 \\
               1 & 0 & 0 & 0 & 0 & 0 & 0 & 0 & 0 \\
             \end{array}\right)$.
\end{center}

\begin{center}
$B_{3}=\left(
         \begin{array}{ccccccccc}
           9 & 6 & 6 & 6 & 4 & 6 & 6 & 6 & 4 \\
           6 & 3 & 4 & 4 & 2 & 4 & 3 & 3 & 2 \\
           6 & 4 & 3 & 3 & 2 & 3 & 4 & 4 & 2 \\
           6 & 4 & 3 & 3 & 2 & 3 & 4 & 4 & 2 \\
           4 & 2 & 2 & 2 & 1 & 2 & 2 & 2 & 1 \\
           6 & 4 & 3 & 3 & 2 & 3 & 4 & 4 & 2 \\
           6 & 3 & 4 & 4 & 2 & 4 & 3 & 3 & 2 \\
           6 & 3 & 4 & 4 & 2 & 4 & 3 & 3 & 2 \\
           4 & 2 & 2 & 2 & 1 & 2 & 2 & 2 & 1 \\
         \end{array}\right)$.
\end{center}

Note that $G_{m,n}$ has the same symmetric transfer matrix in both
horizontal and vertical directions. So the same approach in
\cite{cw} can be taken here. Then we can easily get the following
results. The lower bound of
$\underset{m,n\rightarrow\infty}{\lim}f(m,n)^{1/mn}$ is
$(\frac{\lambda_{p+2q}}{\lambda_{2q}})^{1/p}$ where $\lambda$'s are
the largest eigenvalues of corresponding $T$'s. And the upper bound
of $\underset{m,n\rightarrow\infty}{\lim}f(m,n)^{1/mn}$ is
$(\xi_{2k})^{1/2k}$ where $\xi$'s are the largest eigenvalues of
corresponding $B$'s. Letting $p=1$ , $q=4$ and $k=4$, we have
\begin{center}
$4.631583395...\leq\underset{m,n\rightarrow\infty}{\lim}f(m,n)^{1/mn}\leq5.765456528...$.
\end{center}
Let $f(m,n)$ denote the number of independent sets of the 8.8.4
lattice  graph $G_{m,n}$. Since the number of vertices of $G_{m,n}$
is $4mn - 2m - 2n$, by the similar reason as for the aztec diamonds,
we can see that the entropy constant of the 8.8.4 lattice is
\begin{center}
$(\underset{m,n\rightarrow\infty}{\lim}f(m,n)^{1/mn})^{1/4}$,
\end{center}
which is between $1.467007628...$ and $1.549560101...$.

\textbf{Remark 1.} In this paper we show that for the
 the number of independent sets, the entropy constants of some lattices are the
same as the entropy constants of the corresponding cylindrical and toroidal
lattices. But this phenomenon may disappear for some other models.
As shown in \cite{yyz}, for dimer problem, the entropy constants of quadratic
lattice with cylindrical and toroidal boundaries are different.

\textbf{Remark 2.} To compute the number of independent sets, for many lattices,
using a single transfer matrix is not enough. 
As a consequence, the approach and the concept of
transfer multiplicity introduced in this paper can be used to deal
with more complex lattices.\\

\begin{flushleft}
\textbf{Acknowledgments}
\end{flushleft}
I would like to thank Professors Z. Chen and H. J. Lai for their
helpful comments.

%\vskip 0.3cm \noindent {\bf Acknowledgements} \vskip 0.5cm \noindent

\vskip1cm \noindent
\newcounter{cankao}
\begin{list}
{[\arabic{cankao}]}{\usecounter{cankao}\itemsep=0cm} \centerline{\bf
References}
%\vspace*{0.5cm} \small
\bibitem{b}
{\sc R. J. Baxter}, {\em Exactly Solved Models in Statistical
Mechanics}, (Academic Press, London, 1982).

\bibitem{b2}
{\sc R. J. Baxter}, {\em Planar Lattice Gass with nearest-Neighbor
Exclusion}, Journal Annals of Combinatorics 3(1999) 191-203.

\bibitem{bet}
{\sc R. J. Baxter, I. G. Enting, and S. K. Tsang}, {\em Hard-Square
Lattice Gass}, J. Statist. Phys. 22 465-489, 1980.

\bibitem{cw}
{\sc N. J. Calkin and H. S. Wilf}, {\em The number of independent
sets in a grid graph}, SIAM J. DISCRETE MATH. Vol. 11, No. 1, pp.
54�C60, February 1998.

\bibitem{c}
{\sc T. Y. Chow}, {\em The Q-Spectrum And Spanning Trees of The
Tensor Products of Bipartite Graphs}, Proceedings of the Amer. Math.
Soc., 125(1997), 3155-3161.

\bibitem{mc}
{\sc M. Ciucu}, {\em Perfect matchings, spanning trees, plane
partitions and statistical physics}, Ph.D. thesis, University of
Michigan, Ann Arbor, MI, 1996.

\bibitem{dg}
{\sc C. Domb, M. S. Green (Eds.)}, {\em Phase Transitions and
Critical Phenomena}, Vol. 1, Academic Press, London, 1972.

\bibitem{e}
{\sc K. Engel}, {\em On the Fibonacci number of an m$\times$n
lattice}, Fibonacci Quart., 28 (1990), pp. 72�C78.

\bibitem{f}
{\sc S. R. Finch}, {\em Several constants arising in statistical
mechanics},  Journal Annals of Combinatorics 3(1999) 323-335.

\bibitem{f2}
{\sc S. R. Finch}, {\em Hard Square Entropy Constant}, \S5.12 in
Mathematical Constants. Cambridge, England: Cambridge University
Press, pp. 342-349, 2003.

\bibitem{h} {\sc M. Hall Jr.},
{\em Combinatorial Theory}, John Wiley and Sons, New York, 1986,
p.20 - 23.

\bibitem{k}
{\sc D. E. Knuth}, {\em Aztec diamond, checkerboard graphs, and
spanning trees}, , Journal of Algebraic Combinatorics 6 (3):253-257
1997.

\bibitem{ps}
{\sc P. A. Pearce and K. A. Seaton}, {\em A Classical Theory of Hard
Squares}, J. Statist. Phys. 53, 1061-1072, 1988.

\bibitem{sn}
{\sc S. R. Salinas and J. F. Nagle}, {\em Theory of the phase
transition in the layered hydrogen-bonded SnCl2 2H2 O crystal},
Phys. Rev. B, 9(1974), 4920-4931

\bibitem{v}
{\sc R. S. Varga}, {\em Matrix Iterative Analysis} Prentice-Hall,
Englewood Cliffs, New Jersey, 1962. Theorem 2.1 and Definition 2.2.

\bibitem{w}
{\sc K. Weber}, {\em On the number of stable sets in an m$\times$n
lattice}, Rostock. Math. Kolloq., 34 (1988), pp. 28�C36.

\bibitem{w2}
{\sc H. S. Wilf}, {\em The problem of the kings}, Electron. J.
Combin., 2 (1995), $\sharp$R3, http://www.combinatorics.org.

\bibitem{yyz}
{\sc W. G. Yan, Y.-N. Yeh, and F. J. Zhang}, {\em Dimer problem on the
cylinder and torus}, Physica A, 387(2008), 6069–6078.

\bibitem{zz}
{\sc Z. H. Zhang}, {\em Merrifield-Simmons Index of Generalized Aztec
Diamond and Related Graphs}, MATCH communications in mathematical
and in computer chemistry, 2006, vol.56 no.3.

\end{list}
\end{document}